\documentclass[a4paper]{amsart}
\date{August 02, 2010}
\author[M.~Umehara]{M. Umehara}
\title[Applications of a completeness lemma]{%
      Applications of a completeness lemma \\  in  minimal surface theory \\
      to various classes  of surfaces}
\address[Masaaki Umehara]{%
   Department of Mathematics, Graduate School of Science,\\
   Osaka University,
   Toyonaka, Osaka 560-0043,
   Japan
}
\email{umehara@math.sci.osaka-u.ac.jp}
\author[K.~Yamada]{K. Yamada}
\address[Kotaro Yamada]{%
   Department of Mathematics,
   Tokyo Institute of Technology,
   O-okayama, Meguro, Tokyo 152-8551,
   Japan%
}
\email{kotaro@math.titech.ac.jp}%%\dedicatory{%
\thanks{
  The authors are partially supported by Grant-in-Aid for 
  Scientific Research (A) No.~19204005,
  and (B) No.~21340016, respectively,
  from the Japan Society for the Promotion of Science.
}
\newtheorem{Thm}{Theorem}
\newtheorem*{Thm*}{Theorem}
\newtheorem*{Prop*}{Proposition}
\newtheorem*{completelemma*}{Completeness Lemma}
\newtheorem*{openquestion*}{Question}
\theoremstyle{definition}
\newtheorem{Def}[Thm]{Definition}
\theoremstyle{remark}
\newtheorem{Rem}[Thm]{Remark}
\newenvironment{claim}{%
 \begin{list}{}{%
  \setlength{\leftmargin}{0em}
  \setlength{\rightmargin}{0em}
  \setlength{\itemindent}{\parindent}
  \setlength{\topsep}{0.5\baselineskip}}
  \item[] \it}{\end{list}}

\newcommand{\affiliationone}[1]{\relax}
\newcommand{\affiliationtwo}[1]{\relax}
\renewcommand{\bibname}{\relax}
\newcommand{\claimnote}[1]{#1}
\renewcommand{\phi}{\varphi}
\renewcommand{\epsilon}{\varepsilon}
\newcommand{\op}[1]{{\operatorname{ #1}}}

\newcommand{\R}{\boldsymbol{R}}
\newcommand{\C}{\boldsymbol{C}}
\newcommand{\SL}{\op{SL}}
\newcommand{\SU}{\op{SU}}
\newcommand{\ord}{\op{ord}}
\renewcommand{\Re}{\op{Re}}
\renewcommand{\Im}{\op{Im}}
\newcommand{\trace}{\op{trace}}
\newcommand{\trans}[1]{\vphantom{#1}^t #1}
\newcommand{\cmcone}{CMC-$1$}
\begin{document}
\maketitle

\begin{abstract}
 We give several applications of a lemma on completeness used by 
 Osserman to show the meromorphicity of Weierstrass data for complete
 minimal surfaces with finite total curvature.
 Completeness and weak completeness are defined 
 for several classes of surfaces which admit singular points. 
 The completeness lemma is a useful machinery for 
 the study of completeness in these classes of surfaces.
 In particular, we show that a constant mean curvature one (i.e. CMC-1)
 surface in de Sitter $3$-space
 is complete if and only if it is weakly complete, 
 the singular set is compact and all the ends are conformally
 equivalent to a punctured disk.
\end{abstract}

%--------------------------------------------------------%
% Common Texts                                           %
%--------------------------------------------------------%

\section*{Introduction}
Firstly, we recall the following lemma given in MacLane
\cite{McL} and Osserman \cite{O1} to show properties of complete minimal
surfaces.
We set
\[
     \Delta:=\{z\in \C\,;\,|z|<1\},\qquad\text{and}\qquad
     \Delta^*:=\Delta\setminus \{0\}.
\]
\begin{completelemma*}[{\cite{McL,O1}}, {\cite[Lemma 9.6]{O2}}]
 Let $\omega(z)$ be a holomorphic $1$-form on $\Delta^*$.
 Suppose that the integral
 $
     \int_{\gamma}|\omega|
 $
 diverges to $\infty$ for all paths $\gamma:[0,1)\to \Delta^*$  
 accumulating at the origin $z=0$.
 Then $\omega(z)$ has at most a pole at the origin.
\end{completelemma*}

This lemma plays a crucial role in showing the meromorphicity of
Weierstrass data for complete minimal surfaces with finite total 
curvature.
We show in the first section  that it is also
useful for showing several types of 
completeness of surfaces which may admit
singularities.

In Section~\ref{sec:cor}, we give a further application:
Let $S^3_1(\subset \R^4_1)$ be the de Sitter $3$-space
with metric induced from the Lorentz-Minkowski $4$-space $\R^4_1$, 
which is a simply-connected $3$-dimensional complete Lorentzian manifold with
constant sectional curvature $1$. 
We consider the projection 
\[
     p_L:\SL(2,\C)\longrightarrow 
         S^3_1=\SL(2,\C)/\SU(1,1).
\]
A holomorphic map $F\colon{}\Sigma^2\to\SL(2,\C)$ 
defined on a Riemann surface $\Sigma^2$ is called {\em null\/} if 
$\det (dF/dz)$ vanishes identically on $\Sigma^2$,
where $z$ is an arbitrary complex coordinate on $\Sigma^2$.
For a null holomorphic immersion $F:\Sigma^2\to \SL(2,\C)$,
\[
     f:=p_L\circ F:\Sigma^2\longrightarrow  S^3_1
 \]
gives a spacelike \cmcone{} 
(i.e.\ constant mean curvature one) 
surface with singularities, called a 
{\it \cmcone{} face\/}.
Conversely, a \cmcone{} face $f\colon{}\Sigma^2\to\SL(2,\C)$
is obtained as $f=p_L\circ F$, 
where $F\colon{}\widetilde\Sigma^2\to\SL(2,\C)$ 
is a null holomorphic immersion defined on the universal cover
$\widetilde\Sigma^2$ of $\Sigma^2$
(see \cite{F} and \cite{FRUYY} for details).
The pull-back metric $ds^2_\#$ of the canonical Hermitian metric of 
$\SL(2,\C)$ by the map $F^{-1}$ of taking inverse matrices
gives a single-valued positive definite metric on 
$\Sigma^2$ (see \cite[Remark 1.11]{FRUYY}).
That is, 
\begin{equation}\label{eq:ds-sharp}
 ds^2_{\#}:= \trace (\alpha_{\#})\trans{\overline{(\alpha_{\#})}} =
             \bigl(1+|G|^2\bigr)^2
	     \left|\frac{Q}{dG}\right|^2\qquad 
	     \biggl(
	         \alpha_{\#}:=(F^{-1})^{-1}d(F^{-1})%=-dF\,F^{-1}
             \biggr),
\end{equation}
where $G$ and $Q$ are the hyperbolic Gauss map and the Hopf
differential, 
respectively (cf.\ \cite[(1.17)]{FRUYY}),
and $\trans{(~)}$ denotes the transposition.
A \cmcone{} face $f:\Sigma^2\longrightarrow  S^3_1$
is called {\it weakly complete\/} if $ds^2_\#$ is complete.

On the other hand, we define another completeness as follows:
Let $(N^3,g)$ be a Riemannian or a Lorentzian
manifold in general.
\begin{Def}\label{def:complete} \rm
 We say a $C^\infty$-map $f\colon \Sigma^2 \to N^3$ is
 {\it complete\/} if there exists a symmetric covariant tensor $T$ on
 $\Sigma^2$ with compact support such that $ds^2+T$ gives a complete
 Riemannian metric on  $\Sigma^2$, where $ds^2$ is the pull-back of the 
 metric $g$ on $N^3$, called the 
 {\em first fundamental form}. 
\end{Def}

If $f$ is complete and the singular set is non-empty, 
then the singular set must be compact, by definition.
If $f$ is a \cmcone{} face, completeness implies weak completeness
(cf.\ \cite[Proposition 1.1]{FRUYY}).
This new definition of completeness is a generalization of the
classical one, namely, 
if $f$ has no singular points, then completeness
coincides with the classical notion of completeness in Riemannian
geometry. 
For complete \cmcone{} faces, we have an analogue of the
Osserman inequality for minimal surfaces in $\R^3$ 
(cf.\ \cite[Theorem 0.2]{FRUYY}).
The goal of this paper is to prove the following assertion:
\begin{Thm*}
A \cmcone{} face in $S^3_1$ is complete if and only if it is weakly
 complete, the singular set is compact and every end is
 conformally equivalent to a punctured disk.
\end{Thm*}
The `only if' part has been proved in \cite{FRUYY}.
The proof of the `if' part in this theorem is given in
Section~\ref{sec:cor}.

In this paper, we shall discuss the relationships 
between completeness and weak completeness
on various classes of surfaces with singularities
as applications of the completeness lemma.
The above theorem is the deepest result amongst them.

\section{Applications of the completeness lemma}
\label{sec:complete}
In this section, we shall give several new applications 
of the completeness lemma, 
which show the importance of this kind of
assertion (see Question in Remark \ref{rmk:Q}). 

Let $(N^3,g)$ be a Riemannian $3$-manifold.
It is well-known that the unit tangent bundle $T_1N^3$  of $N^3$
has a canonical contact structure.
A map $L:\Sigma^2\to T_1N^3$
defined on a $2$-manifold $\Sigma^2$
is called a {\it Legendrian immersion\/} if and only if  $L$ is an
isotropic immersion.
Then a map $f:\Sigma^2\to N^3$
is called a {\em wave front\/} or a {\em front\/}
if there exists a Legendrian immersion $L_f:\Sigma^2\to T_1N^3$ 
that is a lift of $f$.
The pull-back metric 
\begin{equation}\label{eq:tau}
   d\tau^2:=L_f^*\tilde g
\end{equation}
of the canonical metric $\tilde g$
of $T_1N^3$  coincides with the sum of the first fundamental form
and the third fundamental form of $f$ if $(N^3,g)$ is a space form.
Then $f$ is called {\it weakly complete\/} if $L_f^*\tilde g$ is 
a complete Riemannian metric on $\Sigma^2$.
It should be remarked that the completeness  
(in the sense of Melko and Sterling \cite{MS})
of surfaces of constant Gaussian curvature $-1$ in $\R^3$
coincides with our notion of weak completeness of wave fronts.
The differential geometry of wave fronts is discussed in
\cite{SUY}.
\begin{Rem}
 This weak completeness is different from that for \cmcone{}
 faces as in the introduction.
 Moreover, weak completeness for improper affine spheres
 given in Remark \ref{rem:impas} of this section is also
 somewhat different from that for fronts and 
 \cmcone{} faces.
 The authors do not yet know of a unified 
 treatment of weak completeness. 
 In fact, there might be several possibilities 
 for completeness of a given
 class of surfaces with singularities.
\end{Rem}

Definition \ref{def:complete} defines completeness for wave fronts.
By definition, completeness implies weak completeness.
These two notions of completeness were defined for flat fronts 
(i.e. wave fronts whose first fundamental form has vanishing Gaussian
curvature on their regular sets) in the hyperbolic $3$-space $H^3$, 
the Euclidean $3$-space $\R^3$, and the $3$-sphere $S^3$,
respectively (cf.\ \cite{KRUY}, \cite{MU} and \cite{KitU}).
In particular, fundamental properties of flat surfaces in 
$H^3$ were given in  G\'alvez, Mart\'\i{}nez and Mil\'an \cite{GMM1}.
Later, further properties for such surfaces as wave fronts were
given in \cite{KUY} and \cite{KRUY}.
 
Let $\Sigma^2$ be a $2$-manifold and $f:\Sigma^2\to H^3$ a flat front.
In this case, there is another lifting of $f$ 
(different from $L_f$) as follows: 
the map $f$ induces a canonical complex structure on $\Sigma^2$, 
and there exists a holomorphic immersion
$F:\widetilde \Sigma^2\to \SL(2,\C)$
such that $f=\pi\circ F$ holds and
$F^{-1}dF$ is off-diagonal,
where 
\[
    \pi:\SL(2,\C)\longrightarrow H^3=\SL(2,\C)/\SU(2)
\]
is the canonical projection.
Remarkably, the weak completeness of $f$ is equivalent to the
completeness of the pull-back metric of the canonical
Hermitian metric of $\SL(2,\C)$ by $F$ (see \cite{KUY}).
For complete flat fronts in $H^3$, we have an analogue of the Osserman
inequality for complete minimal surfaces (see \cite{KUY} for details).
The caustic (i.e. focal surface) 
of a complete flat front $f:\Sigma^2\to H^3$ is weakly
complete (by taking a double cover of $\Sigma^2$ if necessary).
The asymptotic behavior of weakly complete (but not complete) flat
fronts in $H^3$ was analyzed in \cite{KRUY2}.
\begin{Prop*}
 A flat front in $H^3$ is complete if and only if it is weakly complete,
and the singular set is compact and each end of the front is conformally
 equivalent to a punctured disk.
\end{Prop*}

It should be remarked that another characterization of completeness
of $f$ using the total curvature of $d\tau^2$ is given
in \cite[Theorem 3.3]{KRUY}.

\begin{proof}
 The `only if' part has been proved in \cite[Theorem 3.3]{KRUY}, since
 completeness at each end implies that the end is conformally equivalent
 to a punctured disk.
 So we shall prove the converse.
 We use the same notations as in \cite{KRUY}.
 Let $f:\Delta^*\to H^3$ be a flat immersion
 which is weakly complete at
 $z=0$.
 Then, the metric $d\tau^2$ given in \eqref{eq:tau}
 is written as  $d\tau^2=|\omega|^2+|\theta|^2$,
 where $\omega$ and $\theta$ are the
 canonical forms as in \cite[(2.4) and (2.7)]{KRUY}.
 Since $f$ is weakly complete at
 $z=0$, the length $L(\gamma)$  with respect to the metric
 $d\tau^2=|\omega|^2+|\theta|^2$
 diverges to $\infty$ for all paths $\gamma:[0,1)\to \Delta^*$  
 accumulating at the origin. 
 Though $\omega$ is defined only on the universal cover
 $\widetilde\Delta^*$ of $\Delta^*$,
 $|\omega|$ is well-defined on $\Delta^*$.
 Hence there exists $\mu\in[0,1)$ such that
 \begin{equation}\label{eq:omega-exp}
    \omega(z) = z^\mu \hat \omega(z)\,dz,
 \end{equation}
 where $\hat \omega(z)$ is a holomorphic function on $\Delta^*$.
(We have not yet excluded the possibility that 
$\hat \omega(z)$ has an essential singularity at $z=0$.) 
We shall now prove that
 $\hat \omega(z)$ has at most a pole there. For this purpose,
 we set $\rho:=\theta/\omega$.
 Then the singular set is characterized by $\{|\rho|=1\}$
 (cf. \cite[(2.6)]{KRUY}).
 We may assume $|\rho|<1$ on $\Delta^*$ without loss of generality
 by exchanging roles of $\omega$ and $\theta$
 (see \cite[page 270]{KRUY}), if necessary,
 because $f$ is an immersion.
 Then we have that for any path $\gamma$ accumulating at $0$, 
 \[
    \infty=L(\gamma)=
      \int_{\gamma}
      \sqrt{|\omega|^2+|\theta|^2}
      \le 
      \int_{\gamma}
      \sqrt{(1+|\rho|^2)|\omega|^2}
      \leq
      \sqrt{2}\int_{\gamma}|\omega|
      \leq \sqrt{2}
      \int_{\gamma}|\hat \omega(z)|\,|dz|.
 \]
 Then by the completeness lemma, 
 $\hat \omega(z)$ has at most a pole.
 On the other hand, $|\rho|$ is well-defined
 on $\Delta^*$ and $|\rho|<1$ holds on $\Delta^*$.
 Thus there exists
 a real number $\nu\in [0,1)$ such that
 $\rho(z)=z^\nu \hat \rho(z)$, where $\hat \rho(z)$ is a holomorphic function
 on $\Delta^*$ (that is, $\hat \rho(z)$ may have
 essential singularity at the origin).
 We have on $\Delta^*$ that
 \[
       1>|z^{1-\nu}|>|z^{1-\nu}\rho(z)|
             =|z \hat \rho(z)|.
 \]
 Then by the great Picard theorem, $z \hat \rho(z)$ is
 meromorphic at the origin.
 Summing up, $\omega$ and $\theta=\rho\omega$
 are meromorphic $1$-forms on a neighborhood of the origin.
 Then \cite[(3.2)]{KRUY} yields 
 that $d\tau^2$ has finite total curvature, and completeness follows
 from  \cite[Theorem 3.3]{KRUY}.
\end{proof}

\begin{Rem}
 The assumption `conformal equivalency of the end to the unit punctured disc'
 in Proposition cannot be dropped.
 However, in other situations, this condition might not be required.
 In fact, a flat front in $\R^3$ is complete if and only if it is weakly
 complete and the singular set is compact (see \cite[Corollary 4.7]{MU}).
 On the other hand, a  weakly complete immersion may not be complete in
 general.
 In fact, we set
 \[
      ds^2=du^2+2\cos \omega(u,v)\,du\,dv+dv^2
 \]
 on $(\R^2;u,v)$, where $\omega:\R^2\to (0,\pi)$
 is a $C^\infty$-function satisfying
 $
       {\partial^2 \omega(u,v)}/{\partial u \partial v}=0.
 $
 Then it is known that (cf \cite{Kit} or \cite{KitU}) 
 there is a flat immersion
 $f:\R^2\to S^3$ into the unit sphere $S^3$ such that
 $ds^2$ is the first fundamental form and
 $2\sin\omega(u,v)\,du\, dv$ is the second fundamental form.
 So if we set
 \[
      \omega(u,v)
          =\arcsin \left(\frac{e^{-u^2}}2 \right)
             +\arcsin \left(\frac{e^{-v^2}}2 \right),
 \]
 there is an associated flat immersion of $\R^2$ into $S^3$.
 In this case, the length of the curve
 $[0,\infty)\ni t\mapsto (t,-t)\in \R^2$ 
 with respect to $ds^2$ is finite, that is
 \[
     \int_0^\infty \sqrt{2\bigl(1-\cos\omega(t,-t)\bigr)}
     =\int_0^\infty 2\sin \frac{\omega(t,-t)}2 dt
     =\int_0^\infty\exp\left(-t^2\right)dt<\infty.
 \]
 Since the sum of the first and the third fundamental forms is given by
 $du^2+dv^2$, 
 the immersion $f$ is weakly complete, but not complete.
\end{Rem}

\begin{Rem}[\claimnote{Improper affine spheres}]\label{rem:impas}
 Improper affine spheres in $\R^3$  are closely related to flat surfaces
 in $H^3$ (see \cite{Mz2} and also \cite{IM}).
 An improper affine sphere with admissible singular points in the affine
 space $\R^3$  is called an {\em improper affine map\/}
 (see \cite{Mz}).
 Mart\'\i{}nez \cite{Mz} defined completeness 
 for improper affine maps
 and proved an analogue of the Osserman inequality.
 Definition \ref{def:complete} is a generalization of this
 completeness and completeness for wave fronts.
 An improper affine map is called {\it weakly complete} if the metric
 $d\tau^2$ as in \cite[(9)]{Mz} is complete.
 Let $\Sigma^2$ be a Riemann surface
 and $(F,G)$  a pair of holomorphic functions on $\Sigma^2$
 such that
 \begin{enumerate}
  \item \label{item:ia-map-1}  
       $\Re(F\,dG)$ is exact,
  \item \label{item:ia-map-2}
       $|dF|^2+|dG|^2$ is positive definite.
 \end{enumerate}
 Then the induced mapping $f:\Sigma^2\to \R^3=\C\times \R$ 
 given by 
 \[
      f:=\left(
           G+\overline F, 
           \frac{|G|^2-|F|^2}{2}+\Re\left(GF-2\int FdG\right)
         \right)
 \]
 is an improper affine map with $d\tau^2=2(|dF|^2+|dG|^2)$.
 Conversely, any improper affine map 
 is given in this way. The following assertion holds:

\begin{claim}{\it
 An improper affine map in $\R^3$ is a wave front.
 Moreover, it is complete if and only if it is weakly complete, 
 the singular set is compact and all ends  are
 conformally equivalent to a punctured disk.}
\end{claim}

 In fact, $\nu:=\left( \bar F-G, 1\right)$
 gives a Euclidean normal vector field along $f$ (cf. \cite[(8)]{Mz}), 
 and one can easily verify that $f$ is a wave front, that is, 
 $(f,[\nu]):\Sigma^2 \to \R^3\times P^2(\R)$
 gives an immersion if and only if (ii) holds. 
(Nakajo \cite[Proposition 4.3]{N} gave an alternative proof of this fact.)

 To prove the second assertion, 
 we use the same notations as in \cite{Mz}.
 Let $f:\Sigma^2\to \R^3$ be a complete improper affine map.
 By \cite[(9)]{Mz}, we have that
 \begin{align*}
  ds^2&=|dF+dG|^2\le (|dF|+|dG|)^2\\
      &=|dF|^2+|dG|^2+2|dF||dG| 
          \le 2(|dF|^2+|dG|^2)=d\tau^2,
 \end{align*}
 which implies the completeness of $d\tau^2$,
 that is, $f$ is weakly complete.
 On the other hand, by \cite[Proposition 1]{Mz},
 all ends of $f$ are conformally equivalent to a  punctured disk.

 Next, we show the converse.
 Let $f:\Delta^*\to \R^3$ be an improper affine immersion which is weakly
 complete at $z=0$, and for which
 the end $z=0$ is of punctured type.
 Then the data $(F,G)$ corresponding to $f$
 is a pair of holomorphic functions on $\Delta^*$,
 and the length $L(\gamma)$ with respect to the metric
 $d\tau^2=2(|dF|^2+|dG|^2)$ diverges to $\infty$
 for all paths $\gamma:[0,1)\to \Delta^*$  accumulating at the origin.
 Since $f$ is an immersion, $|dF|\ne |dG|$ holds.
 So we may assume
 $|dF|<|dG|$ holds on $\Delta^*$ without loss of generality.
 Then we have that
 \[
   \infty=L(\gamma)=
    \sqrt{2}\int_{\gamma}
      \sqrt{|dF|^2+|dG|^2}
      \le 2\int_{\gamma}|dG|.
 \]
 Thus, by the completeness lemma, $dG(z)$ has at most a pole at the origin.
 Since $|dF/dG|<1$, 
 the great Picard theorem yields that $dF(z)/dG(z)$ has at most a pole at
 the origin.
 In particular, both $dF$ and $dG$ have at most a pole at the
 origin.
 Thus $dF(0)/dG(0)$ is well-defined,
 and satisfies $|dF(0)/G(0)|\le 1$.
 Since  the singular set 
 $\{z\in \Delta^*\,;\, |dF(z)/dG(z)|=1\}$ is empty,
 we have that $|dF(0)/dG(0)|< 1$.
 Then there exists $\delta\in (0,1)$ such that
 $|dF(z)/dG(z)|<\delta$ holds near $z=0$, and
 \begin{align*}
  ds^2&=|dG+dF|^2=|dG|^2\left|1+\frac{dF}{dG}\right|^2
       \ge |dG|^2 \left(1-\left|\frac{dF}{dG}\right|\right)^2 \\
      &\ge (1-\delta)^2 |dG|^2
       \ge \frac{(1-\delta)^2}2\left(|dF|^2+|dG|^2\right),
 \end{align*}
 which proves the completeness of $f$.
\end{Rem}

\section{Proof of Theorem}
\label{sec:cor}
Now we give a proof of the theorem in the introduction.
It is sufficient to show the converse of \cite[Proposition 1.1]{FRUYY}.
We use the same  notations as in \cite{FRUYY}.

Let $f\colon{}\Delta^*\to S^3_1$ be a \cmcone{} immersion
which is weakly complete at the origin.
Since $f$ is an immersion, the metric $d\sigma^2$ given in
\cite[(1.15)]{FRUYY} is a metric of constant curvature $-1$ on $\Delta^*$.
By \cite[Theorem 2.1 and Definition 3.2]{FRUYY}, $f$ is $g$-regular. 
Also, by \cite[Corollary 3.1]{FRUYY}, $f$ must be elliptic or parabolic.
\subsection*{Elliptic case}
In this case, $f$ is a $g$-regular elliptic end.
Since the Schwarzian derivative $S(g)$ of the
secondary Gauss map $g$ is a projective connection
of elliptic type (see \cite[Section 2]{FRUYY}),
$g$ is written in the form 
(see \cite[Proposition 2.2]{FRUYY})
\[
    g = z^{\mu} h(z),
\]
where $\mu$ is a real number and
$h$ is a holomorphic function with $h(0)\neq 0$.
If $\mu\ne 0$, then we can replace $g$ by $1/g$ 
(see \cite[Remark 1.9]{FRUYY}),
and so we may assume that $\mu>0$. In particular, $g(0)=0$ holds.
On the other hand, if $\mu=0$, then
$|g(0)|\ne 1$, since the singular set $\{|g|=1\}$ does not
accumulate at the origin.
Then, 
replacing $g$ by $1/g$ if necessary
(cf.\ \cite[Remark 1.9]{FRUYY}), we may assume
$|g|^2<1-\epsilon$ on $\Delta^*$ for sufficiently small $\epsilon>0$.
When $\mu>0$, the inequality $|g|^2<1-\epsilon$ holds trivially.
Thus we have that
\begin{equation}\label{eq:hat_ds}
    d\hat s^2:=
      (1+|g|^2)^2|\omega|^2 
      \leq 4 |\omega|^2 
      \leq \frac{4}{\epsilon^2} (1-|g|^2)^2|\omega|^2
      =\frac{4}{\epsilon^2}ds^2,
\end{equation}
where $(g,\omega)$ is the Weierstrass data as in \cite[(1.6)]{FRUYY}.

The weak completeness means the completeness of the  metric
$ds^2_\#$ as in \eqref{eq:ds-sharp}.
(The metric $ds^2_\#$ coincides with the metric given in
\cite[(1.17)]{FRUYY}.)
On the other hand, the completeness of $ds^2_\#$ is equivalent to that
of $d\hat s^2$ (see the last part of the proof of 
\cite[Proposition 1.1]{FRUYY}). 
By \eqref{eq:hat_ds}, $ds^2$ is complete at the origin.

\subsection*{Parabolic case}
In this case, by \cite[Lemma P]{FRUYY}, $f$ is a 
$g$-regular
parabolic end of the first kind.
In the proof of Theorem 3.2 in \cite{FRUYY}, 
completeness is not used, but what is applied is the fact that the
ends are immersed, $g$-regular and of punctured type. So we can refer 
to all
of equations in that proof.
One can choose the secondary Gauss map $g$ as
(see \cite[(3.2)]{FRUYY})
\[
   \frac{1}{i}\frac{g(z)+1}{g(z)-1}=\hat g(z)=i(h(z)+ \epsilon \log z),
\]
where $\epsilon=\pm 1$ and $h(z)$ is a holomorphic function on 
a sufficiently small neighborhood of the origin.
Then we have that
\[
   g(z):=\frac{\hat g(z)-i}{\hat g(z)+i},\qquad
   g'(z)=\frac{2(zh'(z)+\epsilon)}{%
           z(h(z)+\epsilon \log z+1)^2},
   	\qquad \left('=\frac{d}{dz}\right).
\]
Since $zh'(z)+\epsilon$ is bounded on a neighborhood
of the origin and
\begin{equation}\label{eq:limit}
    \frac{|(h(z)+\epsilon\log z+1)|}{|\log z|}\to 1  
\end{equation}
as $z\to 0$,
there exists a positive constant  $c_1$
such that
\begin{equation}\label{eq:limit2}
 |g'|\le       \frac{c_1}{|z|\left|\log z\right|^2}.
\end{equation}
On the other hand, by \eqref{eq:limit} we have that
\begin{equation}\label{eq:limit3} 
   |g(z)| = \left|
          \frac{h(z)+\epsilon\log z-1 }{h(z)+\epsilon\log z+1}
         \right|
	  \to 1 \qquad (z\to 0).
\end{equation}
Moreover, since
\[
   1-|g|^2=\frac{4(\Re h+\epsilon \log|z|)}{|\log z+\epsilon(h+1)|^2},
\]
it can be easily checked that
there is a constant $c_2$ such that
\begin{equation}\label{eq:limit4}
  \frac{c_2}{|\log z|}\le |1-|g|^2|.
\end{equation}
Since  $ds^2=(1-|g|^2)^2|Q/dg|^2$, 
\eqref{eq:limit2} yields that there is a constant $c'_2(>0)$  such that
\begin{equation}\label{eq:limit4a}
  |z|^2|\log z|^2
     \left|\frac{Q}{dz}\right|^2 
     \le c'_2 ds^2,
\end{equation}
where $Q=\omega\,dg$ is the Hopf differential
(cf. \cite[(1.8)]{FRUYY}). 
Then, there exist positive constants $c_3$ and $c'_3$ such that
\begin{equation}\label{eq:hats}
    d\hat s^2:=
      (1+|g|^2)^2|\omega|^2 
      \leq c_3 \left|\frac{Q}{dg}\right|^2 
      \leq c'_3 
      \left|\log z\right|^4\left|\frac{Q}{dz}\right|^2.
\end{equation}

We recall the following inequality 
shown on \cite[Appendix A]{FRUYY}
\begin{equation}\label{eq:limit5}
    ds^2 \le c_0 |Q/dz|^2,
\end{equation}
which holds for $g$-regular parabolic ends of the first kind,
where $c_0$ is a positive constant.
By \eqref{eq:limit4}, \eqref{eq:limit4a} \eqref{eq:hats} and
\eqref{eq:limit5},
we have that
\begin{align*}
   d\hat s^2
      &
      \leq c'_3 
           \left|\log z\right|^4\left|\frac{Q}{dz}\right|^2
      \leq  c'_3 c'_2
            \frac{\left|\log z\right|^2}{|z|^2}ds^2
      \leq  c'_3 c'_2 c_0
            \frac{\left|\log z\right|^2}{|z|^2}\left|\frac{Q}{dz}\right|^2\\
     & \leq  
            \frac{c'_3 (c'_2)^2 c_0}{|z|^4}ds^2
      \leq   c'_3 (c'_2 c_0)^2
            \left|\frac{Q}{z^2 dz}\right|^2.
\end{align*}
Again, by the equivalency of completeness of the two metrics $ds^2_\#$ 
and $d\hat s^2$, 
the weak completeness implies the completeness of $d\hat s^2$. 
In particular, 
$|Q/(z^2dz)|^2$ is a complete metric at the origin.
Then by the completeness lemma, $Q/(z^2dz)$ (and so $Q$)
has at most a pole at $z=0$.
We denote by $\ord_0 Q$ the {\em order\/} of $Q$ at the origin.
For example, $\ord_0 Q=m$ if $Q=z^m\,dz^2$.

Firstly, we consider the case that $z=0$ is not a regular end.
Then the hyperbolic Gauss map $G(z)$ of $f$
has an essential singular point
at $z=0$ and the Schwarzian derivative
$S(G)$ has pole of order $<-2$.
Since $f$ is $g$-regular, $S(g)$ is of order $-2$ 
(see \cite[Definition 3.2]{FRUYY}),
so the identity $2Q=S(g)-S(G)$ (cf.\ \cite[(1.20)]{FRUYY}) 
implies $\ord_0 Q\leq -3$.
On the other hand, if $z=0$ is a regular end, 
then $\ord_0 Q= -2$ by \cite[Lemma 5.1]{FRUYY}.
Thus, the Hopf differential of $f$ satisfies 
$\ord_0 Q\le -2$.
In particular, there exists a constant $c_4>0$ such that
\[
     |Q|\ge \frac{c_4|dz|^2}{r^2}\qquad (r:=|z|)
\]
holds on $\Delta^*$.
By \cite[(3.4)]{FRUYY}, we have that
\[
    d\sigma^2:=\frac{4\,|dg|^2}{(1-|g|^2)^2}\le 
                     \frac{C^2|dz|^2}{r^2(c+\log r)^2},
\]
where $C$ and $c$ are positive constants.
Thus there exists a constant $c_5>0$ such that 
(cf.\ \cite[(1.15) and (1.16)]{FRUYY})
\[
 ds^2=4\frac{|Q|^2}{d\sigma^2}
 \ge 4(c_5)^2r^2|Q|^2(c+\log r)^2
 \ge 4(c_4c_5)^2\frac{(c+\log r)^2}{r^2}|dz|^2.
\]
Since $|dz|\geq dr$ and
\[
  \int_1^t \frac{(c+\log r)}{r}dr
       =\frac{(\log t)^2}{2}+c\log t
\]
diverges to $\infty$ as $t\to 0$,
the metric $ds^2$ is complete at $z=0$, 
which proves the assertion.

\begin{Rem}\label{rmk:Q}
Related to the above proof of our main theorem,
we leave here the following:

\begin{openquestion*}
 Let $\omega(z)$ be a holomorphic $1$-form on $\Delta^*$ and
 $n$ a non-negative integer.
 Suppose that the integral
 \[
     \int_{\gamma}
         \left|\omega(z) (\log z)^n \right|
 \]
 diverges to $\infty$ for all  $C^\infty$-paths
 $\gamma:[0,1)\to \Delta^*$  
 accumulating at the origin $z=0$.
 Then does $\omega(z)$ have at most a pole at the origin?
\end{openquestion*}

When $n=0$, this question reduces to the original lemma. 
If the answer is affirmative,
one can obtain the meromorphicity of the Hopf differential
$Q$ directly by applying it to equation
\eqref{eq:hats}, since the weak completeness implies the
completeness of the metric $d\hat s^2$.
The proof of the completeness lemma given in \cite{O2}
cannot be modified directly, since the estimate of 
$\Im(\log z)$ along a path $\gamma$ seems difficult.
Fortunately, in our situation, we have succeeded to prove
our main result without applying the statement, because of \eqref{eq:limit5}.
It should be remarked that this key inequality  \eqref{eq:limit5} 
itself comes from the Gauss equation of the surface.
\end{Rem}

\begin{Rem}[\claimnote{Spacelike maximal surfaces in $\R^3_1$}]
 It is known that \cmcone{} surfaces in $S^3_1$
 have a quite similar properties to spacelike maximal surfaces in 
 the Lorentz-Minkowski space $\R^3_1$ of dimension $3$
 with the metric of signature $(-,+,+)$.
 We consider a fibration
 \[
   p_L:\C^3\ni (\zeta^1,\zeta^2,\zeta^3)\longmapsto 
       \Re \left(-\sqrt{-1} \zeta^3,\zeta^1,\zeta^2\right)\in \R^3_1.
 \]
 The projection of null holomorphic immersions
 into $\R^3_1$ by $p_L$ gives spacelike maximal surfaces with singularities,
 called {\it maxfaces} (see \cite{UY3} for details).
 Here a holomorphic map $F=(F_1,F_2,F_3)\colon\Sigma^2\to \C^3$ defined 
 on a Riemann surface $\Sigma^2$ is called {\em null\/} if 
 $(F_1')^2+(F_2')^2+ (F_3')^2$ vanishes identically, where $'=d/dz$
 denotes the derivative with respect to a local complex coordinate $z$
 of $\Sigma^2$.
 The completeness and the weak completeness for 
 maxfaces are defined in \cite{UY3}.
 As in the case
 of \cmcone{} faces in $S^3_1$, completeness implies
 weak completeness, however, the converse is not true. 
 For complete maxfaces, we have an analogue of the
 Osserman inequality (cf. \cite{UY3}).
 On the other hand, recently, the existence of
 many weakly complete bounded maxfaces in $\R^3_1$
 has been shown, and such surfaces cannot be complete
 (see \cite{MUY}). We can prove the following
assertion by applying the completeness lemma:

\begin{claim}
{\it  A maxface in $\R^3_1$ is complete if and only if it is weakly
  complete, the singular set is compact and all ends  are
  conformally equivalent to a punctured disk.}
\end{claim}

The proof of this assertion 
is easier than the case of CMC-1 faces in $S^3_1$.
 The `only if' part has been proved in \cite[Corollary 4.8]{UY3}, 
 since finiteness of total curvature implies that all ends are
 conformally  equivalent to a punctured disk. 
 So it is sufficient to show the converse.
 The notations are the same as in \cite{UY3}.
 Let $f:\Delta^*\to \R^3_1$ be a spacelike maximal immersion
 which is weakly complete at $z=0$, namely, the length $L(\gamma)$
 with respect to the metric $d\sigma^2=(1+|g|^2)^2|\omega|^2$
 given in \cite[Definition 2.7]{UY3}
 diverges to $\infty$ for all paths $\gamma:[0,1)\to \Delta^*$  
 accumulating at the origin.
 Since $f$ is an immersion, we may assume that $|g(z)|\ne 1$
 for all $z\in \Delta^*$.
 In particular, the image $g(\Delta^*)$ has 
 infinitely many exceptional values.
 Then by the great Picard theorem,
 $g$ has at most a pole at the origin.
 Without loss of generality, we may assume that $g(0)\in \C$.
 Then there exists $M>0$ such that $|g(z)|<M$
 for $z\in \Delta^*$.
 Thus we have that
 \[
   \infty=L(\gamma)=
      \int_{\gamma}d\sigma
      \le (1+M^2)\int_{\gamma}|\omega|.
 \]
 By the completeness lemma, 
 $\omega$ has at most a pole at the origin. 
 Since the Weierstrass data $(g,\omega)$ has
 at most a pole at the origin, $d\sigma^2$ has 
 finite total curvature.
 Then completeness follows from \cite[Corollary 4.8]{UY3}.
\end{Rem}

%%%%%%%%%%%%%%%%%%%%%%%%%%%%%%%%%%%%%%%%%%%%%%%%%%%%%%%%%%%%%%%%%%%%%

\affiliationone{
   M. Umehara\\
   Department of Mathematics,\\
   Graduate School of Science, \\
   Osaka University,\\
   Toyonaka, Osaka 560-0043,\\
   Japan\\
   \email{umehara@math.sci.osaka-u.ac.jp}
}
\affiliationtwo{
   K. Yamada\\
   Department of Mathematics,\\
   Tokyo Institute of Technology,\\
   O-okayama, Meguro, Tokyo 152-8551\\
   Japan
   \email{kotaro@math.titech.ac.jp}
}
\end{document}